\documentstyle{amsppt}
\magnification=\magstep1 \NoBlackBoxes \TagsOnRight
\baselineskip 20pt \pagewidth {5 in} \pageheight {8 in}
\NoRunningHeads \topmatter
\title minimal period estimates for brake orbits of nonlinear
symmetric Hamiltonian systems
\endtitle
\author Chungen Liu \endauthor
\affil School of Mathematics and LPMC, Nankai University\\
Tianjin 300071, People's Republic of China \\
liucg\@nankai.edu.cn\endaffil
\thanks
* Partially supported by the NSF special funds(10531050) and the innovation group funds(10621101),
973 Program of MOST(2006CB805903).
\endthanks
\abstract In this paper,  we consider the minimal period estimates
for brake orbits of nonlinear symmetric Hamiltonian systems.
 We prove
that if the Hamiltonian function $H\in C^2(\Bbb R^{2n}, \Bbb R)$ is
super-quadratic and convex, for every number $\tau>0$, there exists
at least one $\tau$-periodic brake orbit $(\tau,x)$ with minimal
period $\tau$ or $\tau/2$ provided $H(Nx)=H(x)$.
\endabstract
\endtopmatter
\def\<{\langle}
\def\>{\rangle}
\def\r{\Bbb R}
\def\sp{\text{Sp}}
\def\rr{\Bbb R^{2n}}
\def\n{\Bbb N}
\def\z{\Bbb Z}
\def\ga{\gamma}
\def\U{\bold U}
\def\P{\Cal P}
\def\Sp{\text{\rm Sp}}
\def\aa{\alpha}
\def\bb{\beta}
\def\ga{\gamma}

\def\lm{\lambda}

\def\om{\omega}
\def\Om{\Omega}
\def\sg{\sigma}

\def\dm{\diamond}
\def\sm{\setminus}

\def\p{\Cal P(2n)}

\def\c{\Bbb C}

\document

 \head \S1 Introduction and main result\endhead

Let $J=\pmatrix 0&-I\\I&0\endpmatrix$ and $N=\pmatrix
-I&0\\0&I\endpmatrix$ with $I$ being the identity of ${\r}^n$.
Suppose $H\in C^2(\rr,\r)$ satisfying $$ H(Nx)=H(x),\qquad \forall\,
x\in\rr.\tag 1.1$$

We consider the following  problem
$$\cases
&\dot{x}(t) = JH'(x(t)), \\
&x(-t) = Nx(t),  \\
 &x(\tau+t) = x(t), \;\;\forall \,t\in \Bbb R. \endcases\tag 1.2$$

A solution $(\tau,x)$ of (1.2) is a special periodic solution of the
Hamiltonian system in (1.2), we call it a {\it brake orbit} and
$\tau$  the brake period of $x$.

The existence and multiplicity of brake orbits on a given energy
hypersurface was studied by many Mathematicians.  In 1987, P.
Rabinowitz in [26] proved that if $H$ satisfies (1.1), $
\Sigma=H^{-1}(1)$ is star-shaped, and $x\cdot H'(x)\neq 0$ for all
$x\in \Sigma$, then there exists a brake orbit on $\Sigma$. In 1987,
V. Benci and F. Giannoni gave a different proof of the existence of
one brake orbit on $\Sigma$ in [1]. In 1989, A. Szulkin in [27]
proved that there exist at least  $n$ brake orbits on $\Sigma$, if
$H$ satisfies conditions in [26] of Rabinowitz and the energy
hypersurface $\Sigma$ is $\sqrt{2}$-pinched. Long, Zhang and Zhu in
[23] proved that there exist at least $2$ geometrically distinct
brake orbits on any central symmetric strictly convex hypersuface
$\Sigma$. Recently, Z.Zhang and the author of this paper in [16]
proved that there exist at least $[n/2]+1$ geometrically distinct
brake orbits on any  central symmetric strictly convex hypersurface
$\Sigma$, furthermore, there exist at least $n$ geometrically
distinct brake orbits on $\Sigma$ if all brake orbits  on $\Sigma$
are non-degenerate.

 In his pioneering work [24], P. Rabinowitz proposed
a conjecture on whether a superquadratic Hamiltonian system
possesses a periodic solution with a prescribed minimal period. This
conjecture has been deeply studied by many mathematicians. For the
strictly convex case, i.e., $H''(x)>0$, Ekeland and Hofer in [6]
proved that Rabinowtz's conjecture is true. We refer to [3]-[6],[8],
[10], [17]-[19], and reference therein for further survey of the
study on this problem.

  For Rabinowitz' conjecture on the second order Hamiltonian systems, similar results
under various convexity conditions have been proved (cf. [5] and
reference therein). In [17] and [19], under precisely the conditions
 of Rabinowitz, Y. Long proved that for any $\tau>0$ the second order system
$$\ddot x+V'(x)=0$$
 possesses a $\tau$-periodic solution $x$ whose minimal period is
at least $\tau/(n+1)$. Similar result for the first order system
(1.1) is still unknown so far.

It is natural to ask the Rabinowitz's question for the brake orbit
problem: for a superquadratic Hamiltonian function $H$ satisfying
condition (1.1), whether the problem (1.2) possesses a solution
$(\tau,x)$ with prescribed minimal period $\tau$ for any $\tau>0$
(brake orbit minimal periodic problem in short).

In this paper we first consider the  brake orbit minimal periodic
problem for the nonlinear Hamiltonian systems. From Section 3, we
have the following result.
 \proclaim{Theorem 1.1} Suppose  the Hamiltonian function $H$ satisfies the conditions:
\item\quad{(H1)} $H\in C^2(\rr, \r)$ satisfying $H(Nx)=H(x),\;\forall  x\in \rr$.
\item\quad{(H2)} there are constants $\mu>2$  and  $r_0>0$ such that
$$   0<\mu  H(x)\le  H'(x)\cdot x,\quad \forall  |x|\ge r_0.   $$
\item\quad{(H3)} $H(x)=o(|x|^2)$  near $x=0$.
\item\quad{(H4)} $H(x)\ge 0$, $\quad \forall  x\in \rr$.
\item\quad{(H5)} $H''(x)>0$, $\quad \forall  x\in \rr$.

\noindent Then there exists a brake orbit $(\tau,x)$ of problem
(1.2) with minimal brake period $\tau$ or $\tau/2$.
\endproclaim

In fact, in Section 3 a more general theorem is proved (see Theorem
3.1) where the superquadratic condition $(H2)$ is relaxed to
$$H(x)=\frac 12(Bx,x)+\tilde H(x) $$
with $\tilde H$ satisfying condition $(H2)$, and the convexity
condition (H5) is relaxed to

$H''(x(t))\ge 0,\;\;\forall t\in\r$ and
$\int^{\tau/2}_0H''(x(t))\,dt >0$ for all brake orbit $(\tau,x)$.

We also prove some results about the  brake orbit minimal periodic
problems for the second order Hamiltonian systems in Section 3.

 \head \S2 Iteration  inequalities of the
$L_0$-index theory \endhead

We observe that the problem (1.2) can be transformed to the
following Lagrangian boundary value problem
$$\cases &\dot x(t)=JH'(x(t)),\\
& x(0)\in L_0,\quad x(\tau/2)\in L_0,\endcases \tag 2.1$$
 where $L_0=\{0\}\times \r^n\subset \rr$.

An index theory suitable for the  study of problem (2.1) was
established in [13] for any Lagrangian subspace $L$. As usual, we
denote
$$\sp(2n)=\{M\in \Cal {L}(\rr)|\,M^TJM=J\},$$

$$\Cal {P}(2n)=\{\gamma\in C([0,1],\sp(2n))|\,\gamma(0)=I_{2n}\}$$
and
$$\Cal {P}_{\tau}(2n)=\{\gamma\in C([0,\tau],\sp(2n))|\,\gamma(0)=I_{2n}\}.$$
 For a symplectic path
$\gamma\in \Cal P(2n)$, its Maslov-type index associated with a
Lagrangian subspace $L$ is assigned to a pair of integers
$(i_L(\gamma),\nu_L(\gamma))\in\z\times \{0,1,\cdots,n\}$. We call
it the $L$-index of $\gamma$ in short. In [23], the index
$\mu_j(\gamma), \;j=1,2$ was defined for $\gamma\in \Cal P(2n)$, the
$\mu_j$-indices are essentially the special $L$-indices for $L=L_0$
and $L=L_1=\r^n\times\{0\}\subset \rr$ up to a constant $n$,
respectively. In order to estimate the period of a brake orbit, we
need to estimate the $L_0$-index of the iteration path $\gamma^k$
associated to the iterated brake orbit $x^k$.

For reader's convenience, we recall the definition of the
$L_0$-index which was first established in [13]. Some properties for
this index theory  are listed in the appendix below. For
${L_{0}}=\{0\}\oplus \r^n$,
 we define the following two subspaces of $\sp(2n)$ by
 $$ \sp(2n)_{L_{0}}^*=\{M\in \sp(2n)|\, \det V\neq 0\},$$

 $$ \sp(2n)_{L_{0}}^0=\{M\in \sp(2n)|\, \det V= 0\},$$
 for  $M=\pmatrix S & V\\T & U\endpmatrix$.

 Since the space $\sp(2n)$ is path connected, and the $n\times n$
 non-degenerated matrix space has two path connected components,
 one with $\det V>0$, and another
 with $\det V<0$,
 the space $\sp(2n)_{L_{0}}^*$ has two path connected components as
 well. We denote by
 $$\sp(2n)_{L_{0}}^{\pm}=\{M\in \sp(2n)|\, \pm\det V>0 \} $$
 then we have $\sp(2n)_{L_{0}}^*=\sp(2n)_{L_{0}}^+\cup \sp(2n)_{L_{0}}^-$.
  We denote the corresponding symplectic path space by
 $$\p_{L_{0}}^*=\{\gamma\in \p|\, \gamma(1)\in \sp(2n)_{L_{0}}^*\} $$
 and
 $$ \p_{L_{0}}^0=\{\gamma\in \p|\, \gamma(1)\in \sp(2n)_{L_{0}}^0\}. $$
  \proclaim{Definition 2.1} We define the ${L_{0}}$-nullity of any
 symplectic path $\gamma\in \p$ by
 $$\nu_{L_{0}}(\gamma)\equiv\dim\ker_{L_0}(\gamma(1)):=\dim\ker V(1)=n-\text{\rm rank} V(1) $$
 with the $n\times n$ matrix function $V(t)$  defined in (2.1).
 \endproclaim
 We note that  $\text{ rank}\pmatrix
 V(t)\\U(t)\endpmatrix=n$, so the complex matrix
 $U(t)\pm\sqrt{-1}V(t)$ is invertible. We define a complex matrix function by
 $$\Cal Q(t)=[U(t)-\sqrt{-1}V(t)] [U(t)+\sqrt{-1}V(t)]^{-1}. $$
 It is easy to see that the matrix $\Cal Q(t)$ is a unitary matrix for
 any $t\in [0,1]$. We denote by
 $$M_+= \pmatrix 0 & I_n\\ -I_n & 0\endpmatrix, \;\; M_-=\pmatrix 0 & J_n\\
  -J_n & 0\endpmatrix, \;\;J_n=\text{\rm diag}(-1,1,\cdots,1). $$
  It is clear that $M_{\pm}\in \sp(2n)_{L_{0}}^{\pm}$.

  For a path $\gamma\in \p_{L_{0}}^*$, we first adjoin it with a simple symplectic path
  starting from $J=-M_+$, i.e., we define a symplectic path by
  $$\tilde {\gamma}(t)=\cases I\cos \frac{(1-2t)\pi}{2}+J\sin\frac{(1-2t)\pi}{2}, \;\;& t\in [0,1/2];\\
  \gamma(2t-1),\; & t\in [1/2,1].\endcases $$
  then we choose a symplectic path $\beta(t)$ in $\sp(2n)_{L_{0}}^*$ starting from
  $\gamma(1)$ and ending at $M_+$ or $M_-$ according to $\gamma(1)\in\sp(2n)_{L_{0}}^+$
  or $\gamma(1)\in\sp(2n)_{L_{0}}^-$, respectively. We now define a joint path by
  $$\bar{\gamma}(t)=\beta*\tilde \gamma:=\cases \tilde \gamma(2t), \;\;& t\in [0,1/2],\\
  \beta(2t-1),\;\; & t\in [1/2,1].\endcases $$
  By the definition, we see that the symplectic path $\bar\gamma$
  starting from $-M_+$ and ending at either $M_+$ or $M_-$.
  As above, we define
  $$\bar {\Cal Q}(t)=[\bar U(t)-\sqrt{-1}\bar V(t)] [\bar U(t)+\sqrt{-1}\bar V(t)]^{-1}. $$
  for $\bar \gamma(t)=\pmatrix \bar S(t) & \bar V(t)\\\bar T(t) & \bar
  U(t)\endpmatrix$. We can choose a continuous function $\bar
  \Delta(t)$ in $[0,1]$ such that
  $$\det \bar {\Cal Q}(t)=e^{2\sqrt{-1}\bar\Delta(t)}. $$
  By the above arguments, we see that the number $\frac{1}{\pi}(\bar
  \Delta(1)-\bar\Delta(0))\in \z$ and it does not depend on
  the choice of the function $\bar\Delta(t)$.   We note that there
  is a positive continuous function $\rho: [0,1]\to (0, +\infty)$
  such that
  $$\det (\bar{U}(t)-\sqrt{-1}\bar{V}(t))=\rho(t)e^{\sqrt{-1}\bar \Delta(t)}. $$
\proclaim{Definition 2.2} For a symplectic path $\gamma\in
\p_{L_{0}}^*$, we define the ${L_{0}}$-index of $\gamma$ by
$$i_{L_{0}}(\gamma)=\frac{1}{\pi}(\bar
  \Delta(1)-\bar\Delta(0)). \tag 2.2$$
\endproclaim

For a $L_0$-degenerate symplectic path $\gamma\in \p_{L_{0}}^0$, its
${L_{0}}$-index is defined by the infimum of the indices of the
nearby nondegenerate symplectic paths.
 \proclaim{Definition 2.3} For a symplectic path $\gamma\in
\p_{L_{0}}^0$, we define the ${L_{0}}$-index of $\gamma$ by
$$i_{L_{0}}(\gamma)=\sup_{U\in \Cal {N}(\gamma)}\inf_{U}\{i_{L_{0}}(\tilde \gamma)|\,
\tilde \gamma\in U\cap\p_{L_{0}}^*\}, $$ where $\Cal {N}(\gamma)$ is
the set of all open neighborhood of $\gamma$ in $\p_{L_{0}}$.
\endproclaim

 Suppose the continuous symplectic path $\gamma:
[0,1]\to \sp(2n)$
 is the fundamental solution of the following linear Hamiltonian
 system
 $$\dot z(t)=J B(t)z(t) \tag 2.3$$
 with $B(t)$ satisfying $B(t+2)=B(t)$ and $B(1+t)N=NB(1-t))$.
  We define the $L_0$-iteration paths $\gamma^k:[0,k]\to \sp(2n)$ of $\gamma$ by
 $$\gamma^1(t)=\gamma(t), \;t\in [0,1], $$
 $$\gamma^2(t)=\cases & \gamma(t), \;t\in [0,1]\\
 & N\gamma(2-t)\gamma(1)^{-1}N\gamma(1), \;t\in [1,2], \endcases$$
  and in general, for $j\in\n$
  $$\gamma^{2j-1}(t)=\cases & \gamma(t), \;t\in [0,1],\\
 & N\gamma(2-t)\gamma(1)^{-1}N\gamma(1), \;t\in [1,2],\\&\cdots\cdots\\
  & N\gamma(2j-2-t)\gamma(1)^{-1}N\gamma(1)\gamma(2)^{2j-5}, \;t\in [2j-3,2j-2], \\
 & \gamma(t-2j+2)\gamma(2)^{2j-4}, \;t\in [2j-2,2j-1],\endcases \tag 2.4$$

$$\gamma^{2j}(t)=\cases & \gamma(t), \;t\in [0,1],\\
 & N\gamma(2-t)\gamma(1)^{-1}N\gamma(1), \;t\in [1,2],\\&\cdots\cdots\\
  & \gamma(t-2j+2)\gamma(2)^{2j-4}, \;t\in [2j-2,2j-1], \\
 & N\gamma(2j-t)\gamma(1)^{-1}N\gamma(1)\gamma(2)^{2j-3}, \;t\in [2j-1,2j],\endcases \tag 2.5$$
We note that if $\tilde \gamma(t),\;t\in \r$ is the fundamental
solution of the linear system (2.3), then there holds
$\gamma^k=\tilde\gamma|_{[0,k]}$.
 For the iteration path $\gamma^k$, the following
Bott-type iteration formulas were proved in [16].

 \proclaim{Proposition 2.4} Suppose $\omega_k=e^{\pi \sqrt{-1}/k}$.
For odd $k$ we have
$$\aligned & i_{L_0}(\gamma^k)=i_{L_0}(\gamma^1)+\sum_{i=1}^{\left[k/2\right]}i_{\omega_k^{2i}}(\gamma^2),
\\ &
\nu_{L_0}(\gamma^k)=\nu_{L_0}(\gamma^1)+\sum_{i=1}^{\left[k/2\right]}\nu_{\omega_k^{2i}}(\gamma^2),\endaligned$$
 for even $k$, we have
$$ \aligned & i_{L_0}(\gamma^k)=i_{L_0}(\gamma^1)+i^{L_0}_{\omega_k^{k/2}}(\gamma^1)+
\sum_{i=1}^{k/2-1}i_{\omega_k^{2i}}(\gamma^2),\;
\\
&
\nu_{L_0}(\gamma^k)=\nu_{L_0}(\gamma^1)+\nu^{L_0}_{\omega_k^{k/2}}(\gamma^1)+
\sum_{i=1}^{k/2-1}\nu_{\omega_k^{2i}}(\gamma^2),
\endaligned$$
where the $(L_0$-$\omega)$ index
$(i_{\omega}^{L_0}(\gamma),\nu_{\omega}^{L_0}(\gamma))$ of $\gamma$
for $\omega\in {\bold U}:=\{z\in\c|\,|z|=1\}$ was defined in [16],
and the $\omega$-index
$(i_{\omega}(\gamma^2),\nu_{\omega}(\gamma^2))$ of $\gamma^2$ for
$\omega\in \bold U$ was defined in [22](cf.[21]).
\endproclaim

We note that $\omega_k^{k/2}=\sqrt{-1}$.
 For any  two $2k_i\times 2k_i$ matrices
of square block form, $M_i=\pmatrix A_i & B_i\\
                                            C_i & D_i \endpmatrix$ with $i=1, 2$,
                                             the $\dm$-product of
$M_1$ and $M_2$ is defined  to be the $2(k_1+k_2)\times 2(k_1+k_2)$
matrix
$$ M_1\dm M_2=\pmatrix A_1 & 0 & B_1 & 0\\
                                               0 & A_2 & 0 & B_2\\
                                           C_1 & 0 & D_1 & 0\\
                                                0 & C_2 & 0 & D_2 \endpmatrix.  $$
Denote by $M^{\dm k}$ the $k$-fold $\dm$-product of $M$. Let
$N_1(\lm, b)=\pmatrix \lm & b\\
                                       0 & \lm\endpmatrix$ for $\lm=\pm 1$ and $b=\pm 1$, or
                                       $0$.
Denote by $R(\theta)=\pmatrix \cos \theta & -\sin \theta\\
                                       \sin\theta & \cos\theta\endpmatrix$. We remind that the unit circle in the
 complex plane is defined by $\bold
U=\{z\in \Bbb C:|z|=1\}$, and the upper(lower) semi closed  unit
circle $\bold U^+$($\bold U^-$) is defined by $\bold U^{\pm}=\{z\in
\bold U: z=e^{\theta\sqrt{-1}},\;0\le \pm \theta\le \pi\}$. In [22],
for any $M\in \Sp(2n)$, Long defined the {\it homotopy set} of $M$
in $\Sp(2n)$ by
$$\aligned
  \Om(M)=\{N\in \Sp(2n)\,|\,&\sg(N)\cap\U=\sg(M)\cap U \; \text{and} \\
   & \dim_{\c}\ker_{\c}(N-\lm I)=\dim_{\c}\ker_{\c}(M-\lm I),\; \forall \lm\in \sg(M)\cap\U\}.\endaligned $$
The path connected component of $\Om(M)$ which contains $M$ is
denoted by $\Om_0(M)$, and is called the {\it homotopy component} of
$M$ in $\Sp(2n)$.

In [15], the following result was proved (cf. Theorem 10.1.1 of
[21]).

 \proclaim{Proposition 2.5} $1^o$ For any $\gamma\in\P(2n)$ and $\om\in{\U}\sm\{1\}$,
there always holds $$ i_{1}(\ga) + \nu_{1}(\ga) - n \le i_{\om}(\ga)
           \le i_{1}(\ga) + n - \nu_{\om}(\ga).       \tag 2.6$$

  $2^o$ The left equality in (2.6) holds for some $\om\in\U^+\sm\{1\}$ (or $\U^-\sm\{1\}$)
if and only if there holds $I_{2p}\dm N_1(1,-1)^{\dm q}\dm
K\in\Om_0(\ga(\tau))$ for some non-negative integers $p$ and $q$
satisfying $0\leq p+q\le n$ and $K\in\Sp(2(n-p-q))$ with
$\sg(K)\subset\U\sm\{1\}$ satisfying that all eigenvalues of $K$
located within the  arc between $1$ and $\om$ including $\om$ in
$\U^+$ (or $\U^-$) possess total multiplicity $n-p-q$. If $\om\ne
-1$, all eigenvalues of $K$ are in $\U\sm \r$ and those in $\U^+\sm
\r$ (or $\U^-\sm \r$) are all Krein negative (or positive) definite.
If $\om =-1$, it holds that $-I_{2s}\dm N_1(-1,1)^{\dm t}\dm H\in
\Om_0(K)$ for some non-negative integers $s$ and $t$ satisfying
$0\le s+t\le n-p-q$, and some $H\in \Sp(2(n-p-q-s-t))$ satisfying
$\sg (H)\subset \U\sm \r$ and that all elements in $\sg (H)\cap
\U^+$ (or $\sg (H)\cap \U^-$) are all Krein-negative (or
Krein-positive) definite.

$3^o$ The left equality in (2.6) holds for all $\om\in\U\sm\{1\}$ if
and only if $I_{2p}\dm N_1(1, -1)^{\dm (n-p)}\in\Om_0(\ga(\tau))$
for some integer $p\in [0,n]$. Specifically in this case, all the
eigenvalues of $\ga(\tau)$ equal to $1$ and $\nu_{\tau}(\ga)=n+p\ge
n$.

$4^o$ The right equality in (2.6) holds for some
$\om\in\U^+\sm\{1\}$ (or $\U^-\sm\{1\}$) if and only only if there
holds $I_{2p}\dm N_1(1,1)^{\dm r}\dm K\in\Om_0(\ga(\tau))$ for some
non-negative integers $p$ and $r$ satisfying $0\leq p+r\le n$ and
$K\in\Sp(2(n-p-r))$ with $\sg(K)\subset\U\sm\{1\}$ satisfying the
condition that all eigenvalues of $K$ located within the closed arc
between $1$ and $\om$ in $\U^+\sm\{1\}$ (or $\U^-\sm\{1\}$) possess
total multiplicity $n-p-r$. If $\om\not=-1$, all eigenvalues in
$\sg(K)\cap\U^+$ (or $\sg(K)\cap\U^-$) are all Krein positive (or
negative) definite; if $\om=-1$, there holds $(-I_{2s})\dm
N_1(-1,1)^{\dm t}\dm H \in \Om_0(K)$ for some non-negative integers
$s$ and $t$ satisfying $0\le s+t\le n-p-r$, and some
$H\in\Sp(2(n-p-r-s-t))$ satisfying $\sg(H)\subset \U\sm\r$ and that
all elements in $\sg(H)\cap\U^+$ (or $\sg(H)\cap\U^-$) are all Krein
positive (or negative) definite.

$5^o$ The right equality in (2.6) holds for all $\om\in \U\sm\{1\}$
if and only if $I_{2p}\dm N_1(1, 1)^{\dm (n-p)}\in\Om_0(\ga(\tau))$
for some integer $p\in [0,n]$. Specifically in this case, all the
eigenvalues of $\ga(\tau)$ must be $1$, and there holds
$\nu_{\tau}(\ga)=n+p\ge n$.

$6^o$ Both equalities in (2.6) hold for all $\om\in \U\sm\{1\}$ if
and only if $\ga(\tau)=I_{2n}$.
\endproclaim

Combining Propositions 2.4 and 2.5, we have the following result.
 \proclaim{Theorem 2.6} $1^o$ For any $\ga\in \Cal P(2n)$ and $k\in\n$,  there holds
$$\align
   & i_{L_0}(\gamma^1)+\left[\frac k2\right] (i_{1}(\ga^2) + \nu_{1}(\ga^2) - n)  \le i_{L_0}(\ga^k) \\
       &\le i_{L_0}(\gamma^1)+\left[\frac k2 \right](i_{1}(\ga) + n )- \frac 12 \nu_{1}(\ga^{2k})+\frac
       12\nu_1(\ga^2), \;\text{if}\;k\in 2\n-1, \tag 2.7\\
   & i_{L_0}(\gamma^1)+i^{L_0}_{\sqrt{-1}}(\gamma^1)+\left(\frac k2-1\right) (i_{1}(\ga^2)
    + \nu_{1}(\ga^2) - n)  \le i_{L_0}(\ga^k) \le i_{L_0}(\gamma^1)+i^{L_0}_{\sqrt{-1}}(\gamma^1)\\
       &+\left(\frac k2-1 \right)(i_{1}(\ga) + n )-
        \frac 12 \nu_{1}(\ga^{2k})+\frac
       12\nu_1(\ga^2)+\frac
       12\nu_{-1}(\ga^2), \;\text{if}\;k\in 2\n. \tag 2.8
\endalign $$
The index $(i_{\omega}^{L_0}(\gamma),\nu_{\omega}^{L_0}(\gamma))$ is
defined in [16] for $\omega\in \Bbb U=\{z\in \Bbb C|\;|z|=1\}$, see
also Definition 4.9 in the appendix below.

$2^o$ The left equality of (2.7) holds for some $k\ge 3$ and of
(2.8) holds for some $k\ge 4$ if and only if there holds $I_{2p}\dm
N_1(1,-1)^{\dm q}\dm K\in\Om_0(\ga^2(2))$ for some non-negative
integers $p$ and $q$ satisfying $p+q\le n$ and some
$K\in\Sp(2(n-p-q))$ satisfying $\sg(K)\subset\U\sm\r$. If
$r=n-p-q>0$, then $R(\theta_1)\dm \cdots \dm R(\theta_r)\in
\Omega_0(K)$ for some $\theta_j\in (0,\pi)$. In this case, all
eigenvalues of $K$ on $\U^+$ (on $\U^-$) are located on the  arc
between $1$ and $\exp(2\pi\sqrt{-1}/k)$ (and
$\exp(-2\pi\sqrt{-1}/k)$) in $\U^+$ (in $\U^-$) and are all Krein
negative (positive) definite.

$3^o$ The right equality of (2.7) holds for some $k\ge 3$ and of
(2.8) holds for some $k\ge 4$ if and only if there holds $I_{2p}\dm
N_1(1,1)^{\dm r}\in\Om_0(\ga^2(2))$ for some non-negative integers
$p$ and $r$ satisfying $p+r= n$.

$4^o$ Both  equalities of (2.7), and also of (2.8), hold for some
$k>2$ if and only if $\ga^2(2)=I_{2n}$.
\endproclaim
\demo{Proof} By Proposition 2.4, summing the inequalities of (2.6)
with $\omega=\omega_k^{2i}, 1\le i<k/2, i\in \n$, we obtain the
inequalities (2.7) for odd $k$ and (2.8) for even $k$. We remind
that here we have used the Bott-type formula
$$\nu_1(\gamma^{k})=\sum_{\omega^k=1}\nu_{\omega}(\gamma).$$
The equality conditions follow from $2^o$ and $4^o$ of Proposition
2.5 together with Corollary 9.2.8 and List 12  in P198 of [21]. We
note that from List 12 in P198 of [21], no eigenvalue on $\bold U^+$
is Krein positive definite.
\enddemo

 Since we should consider the Bott-type iteration formulas in
Proposition 2.4 in  odd and even cases, the inequalities in Theorem
2.6 is naturally considered in two cases correspondingly.  We will
see that the inequalities in Theorem 2.6 for even times iteration
path are our main difficult to prove that the brake orbit found in
Section 3 has minimal period, though we believe  this kind brake
orbit has minimal period, we can only prove that it has minimal
period or it is $2$-times iteration of a brake orbit with minimal
period.

\head \S3 Applications to nonlinear Hamiltonian systems \endhead

We now apply Theorem 2.6 to the brake orbit problem of  autonomous
Hamiltonian system
$$ \cases &-J\dot x=Bx+H'(x), \qquad x\in \rr,\\
& x(\tau/2+t)=Nx(\tau/2-t),\\
& x(\tau+t)=x(t), \;t\in\Bbb R,\endcases \tag 3.1 $$
where $H(Nx)=H(x)$  and $B=\pmatrix B_1 & 0 \\
0 & B_2\endpmatrix$ is a $2n\times 2n$ symmetric semi-positive
definite matrix whose operator norm is denoted by $\|B\|$, $B_1$ and
$B_2$ are $n\times n$ symmetric matrices. A solution $(\tau,x)$ of
the problem (3.1) is a  brake orbit of the Hamiltonian system, and
$\tau$ is  the  brake period of $x$. To find a brake orbit of the
Hamiltonian system in (3.1), it is sufficient to solve
 the following problem
$$\cases &-J\dot x(t)=Bx+H'(x(t)), \;\; x\in \rr,\;t\in [0,\tau/2],\\
& x(0)\in L_0,\; x(\tau/2)\in L_0. \endcases \tag 3.2  $$
 Any solution $x$ of problem (3.2) can be extended to a brake orbit
 $(\tau,x)$ with the mirror symmetry of $L_0$ by $x(\tau/2+t)=Nx(\tau/2-t),\;t\in [0,\tau/2]$
 and $x(\tau+t)=x(t),\;t\in\r$.
\proclaim{Theorem 3.1} Suppose  the Hamiltonian function $H$
satisfies the conditions:
\item\quad{(H1)} $H\in C^2(\rr, \r)$ satisfying $H(Nx)=H(x),\;\forall  x\in \rr$.
\item\quad{(H2)} there are constants $\mu>2$  and  $r_0>0$ such that
$$   0<\mu  H(x)\le  H'(x)\cdot x,\quad \forall  |x|\ge r_0.   $$
\item\quad{(H3)} $H(x)=o(|x|^2)$  at $x=0$.
\item\quad{(H4)} $H(x)\ge 0$ $\quad \forall  x\in \rr$.

\noindent Then for every $0<\tau<\frac{2\pi}{\|B\|}$, the system
(3.1) possesses a non-constant brake orbit $(\tau,x)$ satisfying
$$ i_{L_0}(x, \tau/2)\le 1.  \tag 3.3$$
\noindent Moreover, if $x$ further satisfies the following
condition:
\item\quad{(HX)} $H''(x(t))\ge 0 \;\;\forall t\in\r$ and $\int^{\tau/2}_0H''(x(t))\,dt >0$.

\noindent Then  the minimal brake period of  $x$  is $\tau$ or
$\tau/2$.
\endproclaim
We remind that if $B=0$, then $\frac{2\pi}{\|B\|}=+\infty$.

 \demo{Proof} We divide the proof into two steps.

 \noindent{\it Step 1.} Show that there exists a brake orbit $(\tau,x)$
 satisfying (3.3) for $0<\tau<\frac{2\pi}{\|B\|}$.

 Fix $\tau\in (0,\frac{2\pi}{\|B\|})$. Without loss
generality, we suppose $\tau=2$, then $\tau<\frac{2\pi}{\|B\|}$
implies $\|B\|<\pi$.   By conditions (H1)-(H4), we can find a
non-constant $\tau$-periodic solution $x$ of (3.2) via the saddle
point theorem such that (3.3) holds. For reader's convenience, we
sketch the proof here and refer the reader to Theorem 3.5 of [15]
for the case of periodic solution. We note that the main ideas here
are the same as that in the periodic case. We refer the paper [11]
for some details.

  In fact, following P. Rabinowitz' pioneering work [24], let $K>0$ and $\chi\in C^{\infty}(\r, \r)$
such that $\chi(t)=1$ if $t\le K$, $\chi(t)=0$ if $t\ge K+1$, and
$\chi'(t)<0$ if $y\in (K, K+1)$. The number $K$ will be determined
later. Set
$$ \hat H_K(z)=\frac{1}{2}(Bz, z) + H_K(z),  $$
with
$$ H_K(z)=\chi(|z|) H(z)+(1-\chi(|z|))R_K|z|^4, $$
where the constant $R_K$ satisfies
$$ R_K\ge \max_{K\le |z|\le  K+1}\frac {H(z)}{|z|^4}. $$

  We set $L^2=L^2([0,1],\rr)$ and define a Hilbert space $E:=\Cal
W_{L_0}=W^{1/2,2}_{L_0}([0,1],\rr)$ with $L_0$ boundary conditions
by
$$\Cal  W_{L_0}=\{z\in L^2|\; z(t)=\sum_{k\in\z}\exp(k\pi tJ)a_k,\; a_k\in L_0,\;\|z\|^2:=\sum_{k\in\z}(1+|k|)|a_k|^2<\infty\}. $$
We denote its inner product by $\<\cdot,\cdot\>$. By the well-known
Sobolev embedding theorem, for any $s\in [1,+\infty)$, there is a
constant $C_s>0$ such that
$$\|z\|_{L^s}\le C_s \|z\|, \;\;\forall \,z\in \Cal  W_{L_0}. $$ Define a
functional $f_K$ on $E$ by
$$ f_K(z)=\int^{1}_0(\frac {1}{2}\dot z\cdot Jz - \hat{H}_K(z))\,dt, \;\;\forall z\in E. \tag 3.4$$
For $m\in\Bbb N$, define $E^0=L_0$,
$$\aligned
   E_m &= \{\,z\in E\;|\; z(t)=\sum_{k=-m}^{m}\exp(k\pi tJ)a_k,\;a_k\in L_0\}, \\
   E^{\pm} &= \{\,z\in E\;|\; z(t)=\sum_{\pm k>0} \exp(k\pi tJ)a_k,\;a_k\in L_0\}, \\
\endaligned$$
and $E^+_m=E_m\cap E^+$,  $E^-_m=E_m\cap E^-$.  We have
$E_m=E^-_m\oplus E^0\oplus E^+_m$.  Let $P_m$ be the projection
$P_m: E\to E_m$. Then $\{E_m,P_m\}_{m\in\Bbb N}$ form a Galerkin
approximation scheme of the operator $-Jd/dt$ on $E$. Denote by
$f_{K,m}=f_K|_{E_m}$. Set $Q_m=\{re:0\le r\le r_1\}\oplus
\{B_{r_1}(0)\cap (E_m^-\oplus E_m^0)\}$ with some $e\in \partial
B_1(0)\cap E_m^+$.  Then  for large $r_1>0$ and small $\rho>0$,
$\partial Q_m$ and $B_{\rho}(0)\cap E^+_m$ form a topological (in
fact homologically) link (cf. P84 of [2]). By the condition
$\|B\|<\pi$, we obtain a constant $\beta=\beta(K)>0$ such that

\item {(I)} $\quad   f_{K,m}(z)\ge \beta>0, \qquad \forall
z\in
\partial B_{\rho}(0)\cap E^+_m$,

\item {(II)} $\quad   f_{K,m}(z)\le 0,\qquad  \forall z\in
\partial Q_m$.

In fact, by (H3), for any $\varepsilon> 0$, there is a $\delta> 0$
such that ${H}_{K}(z)\leq \varepsilon|z|^{2}$ if $|z|\leq \delta$.
Since $\hat{H}_{K}(z)|z|^{-4}$ is uniformly bounded as
$|z|\rightarrow +\infty$, there is an $M_{1}=M_{1}(K)$ such that
$\hat{H}_{K}(z)\leq M_{1}|z|^{4}$ for $|z|\geq \delta$. Hence
$$
\hat{H}_{K}(z)\leq \varepsilon|z|^{2}+M_{1}|z|^{4},\quad \forall
z\in \rr.
$$
For $z \in \partial B_{\rho}(0)\cap E^+_{m}$, we have
$$
\int_{0}^{1} {H}_{K}(t,z)dt\leq
\varepsilon\|z\|^{2}_{L^{2}}+M_{1}\|z\|^{4}_{L^{4}}\leq (\varepsilon
C_{2}^{2}+M_{1}C_{4}^{4}\|z\|^{2})\|z\|^{2}.
$$
So we have
$$ \aligned f_{K,m}(z) &= \frac{1}{2} \langle Az,z \rangle
-\frac{1}{2}\langle {B}z, z\rangle-\int_{0}^{1} {H}_{K}(z(t))dt \\
 &\geq
 \frac{\pi}{2}\|z\|^{2}-\frac{\|B\|}{2}\|z\|^{2}-
(\varepsilon C_{2}^{2}+M_{1}C_{4}^{4}\|z\|^{2})\|z\|^{2} \\
&= \frac{\pi}{2}\rho^{2}-\frac{\|B\|}{2}\rho^{2}-(\varepsilon
C_{2}^{2}+M_{1}C_{4}^{4}\rho^{2})\rho^{2}. \endaligned$$
 Since
$\|B\|<\pi$, we can choose constants $\rho=\rho(K)>0$ and
$\beta=\beta(K)>0$, which are sufficiently small and independent of
$m$, such that for $z\in \partial B_{\rho}(0)\cap E^+_{m}$,
$$
f_{K,m}(z)\geq \beta>0.
$$
Hence (I) holds.

 Let $e\in E_{m}^{+}\cap \partial B_{1}$ and
$z=z^{-}+z^{0}\in E^-_m\oplus E^0$. We have
$$\aligned
f_{K,m}(z+re) &= \frac{1}{2}\langle Az^{-}, z^{-}
\rangle+\frac{1}{2}r^{2}\langle
Ae, e \rangle-\frac12 \langle {B}(z+re), z+re\rangle-\int_{0}^{1} \hat{H}_{K}(z+re)dt \\
&\leq -\frac{\pi}{2}\|z^{-}\|^{2}+
\frac{\pi}{2}r^{2}-\int_{0}^{1}\hat{H}_{K}(z+re)dt,
\endaligned$$
If $r=0$, from condition (H4), there holds
$$f_{K,m}(z+re)\le -\frac{\pi}{2}\|z^{-}\|^{2}\le 0.$$
If $r=r_1$ or $\|z\|=r_1$, then from (H2), We have
$$
{H}_{K}(z)\geq b_{1}|z|^{\mu}-b_{2},
$$
 where
$b_{1}>0,\,b_{2}$ are two constants independent of $K$ and $m$. Then
there holds
$$\aligned
\int_{0}^{1} \hat{H}_{K}(z+re)dt &\geq b_{1}\int_{0}^{1}
|z+re|^{\mu}dt-b_{2}  \\ &\geq b_{3}\left(\int_{0}^{1}
|z+re|^{2}dt\right)^{\frac{\mu}{2}}-b_{4}  \\
&\geq b_{5}\left(\|z^{0}\|^{\mu}+r^{\mu}\right)-b_{4},
\endaligned$$
where $b_{3},\,b_{4}$ are constants and $b_{5}>0$ independent of $K$
and $m$. Thus there holds
$$
f_{K,m}(z+re)\leq -\frac{\pi}{2}\|z^{-}\|^{2}+ \frac{\pi}{2}r^{2}-
b_{5}\left(\|z^{0}\|^{\mu}+r^{\mu}\right)+b_{4},
$$
So we can choose large enough $r_{1}$  independent of $K$ and $m$
such that
$$
\varphi_{m}(z+re)\leq 0, \quad {\text on }\; \partial Q_{m}.
$$
Then (II) holds.

Now define $\Omega=\{\Phi\in C(Q_m, E_m)\,|\,\Phi(x)=x\;\text{
for}\;x\in{\partial Q_m}\}$, and set
$$c_{K,m}=\inf_{\Phi\in\Omega}\sup_{x\in \Phi(Q_m)}f_{K,m}(x).$$

It is well known that $f_K$ satisfies the usual (P.S)$^{\ast}$
condition on $E$, i.e. a sequence $\{x_m\}$ with $x_m\in E_m$
possesses a convergent subsequence in $E$, provided
$f_{K,m}'(x_m)\to 0$ as $m\to\infty$ and $|f_{K,m}(x_m)|\le b$ for
some $b>0$ and all $m\in\n$ (see [11] for a proof). Thus by the
saddle point theorem (cf. [25]), we see that $c_{K,m}\ge \beta>0$ is
a critical value of $f_{K,m}$, we denote the corresponding critical
point by $x_{K,m}$. The Morse index of $x_{K,m}$ satisfies
$$m^-(x_{K,m})\le \dim Q_m=mn+n+1.$$
By taking $m\to +\infty$, we obtain a critical point $x_K$ such that
$x_{K,m}\to x_K$, $m\to +\infty$ and $m_d^-(x_K)\le \dim
Q_m=1+n+mn$,
  $0< c_K\equiv f_K(x_K)\le M_1$, where the $d$-Morse index $m^-_{d}(x_K)$ is defined to the total number
  of the eigenvalues of $f''_{K}$ belonging to $(-\infty,d]$ for $d>0$ small enough, and
  $M_1$ is a constant
independent of $K$. Moreover,  by the Galerkin approximation method,
Theorem 2.1 of [14], we have the $d$-Morse index satisfying

$$m_d^-(x_K)=mn+n+i_{L_0}(x_K,1)\le 1+n+mn.$$
Thus we have
$$i_{L_0}(x_K,1)\le 1. $$

    Now the similar arguments as in the section 6 of [25] yields a constant $M_2$ independent of $K$
such that $\|x_K\|_{\infty}\leq M_2$. Choose $K>M_2$. Then $x\equiv
x_K$ is a non-constant solution of the problem (3.2) satisfying
(3.3). By extending the domain with mirror symmetry of $L_0$, we
obtain a $2$-periodic brake orbit $(2,x)$ of problem (3.1).

 {\it Step 2.}  Estimate the brake period of $(2, x)$.

  Denote the minimal period of the brake orbit $x$ by $2/k$ for some $k\in \Bbb
   N$, i.e., $(x, 1/k)$ is a solution of the problem (3.2).
By the condition (HX) and $B$ being semi-positive definite, using
(9.17) of [4],  we have that $i_1(x,2/k)\ge n$ for every
$2/k$-periodic solution $(x,2/k)$ (see also (4.2) in the appendix
below), and by Theorem 5.2 of [13] we see that $i_{L_0}(x,1/k)\ge 0$
for the $L_0$-solution $(x,1/k)$ (see also (4.3) in the appendix
below). Together with (3.21) of [16](see also (4.14) in the appendix
below), we obtain
$$i^{L_0}_{\sqrt{-1}}(x,1/k)\ge i_{L_0}(x,1/k)\ge 0,\;\;i_1(x,2/k)\ge n. \tag 3.5$$
Since the system (3.1) is autonomous, we have
$$ \nu_{1}(x,2/k)\ge 1. \tag 3.6$$
Therefore, by Theorem 2.6, (3.3) and (3.5)-(3.6), we obtain $k=1,
2,3,4$.

If $k=3$, by (3.3) and (3.5)-(3.6), and by using Theorem 2.6 again
we find the left
 equality of (2.7) holds for $k=3$ and $i_{L_0}(x,1/3)= 0,
i_1(x,2/3)= n$, and $\nu_1(x,2/3)= 1$.

 The left side hand
equality in the  inequality (2.7) holds if and only if $I_{2p}\dm
N_1(1,-1)^{\dm q}\dm K\in\Om_0(\ga(2/3))$ for some non-negative
integers $p$ and $q$ satisfying $p+q\le n$ and some
$K\in\Sp(2(n-p-q))$ satisfying $\sg(K)\subset\U\sm\r$. If
$r=n-p-q>0$, then by List 12 in P198 of [21](see also the list after
Definition 4.4 in the appendix below), we have $R(\theta_1)\dm
\cdots \dm R(\theta_r)\in \Omega_0(K)$ for some $\theta_j\in
(0,\pi)$. In this case, all eigenvalues of $K$ on $\U^+$ (on $\U^-$)
are located on the  arc between $1$ and $\exp(2\pi\sqrt{-1}/k)$ (and
$\exp(-2\pi\sqrt{-1}/k)$) on $\U^+$ (in $\U^-$) and are all Krein
negative (positive) definite. We remind that $\gamma(t)$ is the
fundamental solution of the linearized system at $(2/3,x)$. By the
condition $\nu_1(x,2/3)=1$, we have $p=0, q=1$. By Lemma 4.3 in the
appendix below, there are paths $\aa\in\P_{2/3}(2)$,
$\bb\in\P_{2/3}(2n-2)$ such that $\ga\thicksim\aa\dm\bb$,
$\aa(2/3)=N_1(1,-1)$, $\bb(\tau)=K$. By the locations of the end
point matrix $\aa(2/3)$ and $\bb(2/3)$, there are two integers $k_1,
\;k_2$ such that (see the proof of Theorem 4.3 in [15], specially
(4.18) and (4.19) there).
$$i_1(\aa,2/3)=2k_1,\;i_1(\bb, 2/3)=2k_2+n-1.$$
From this result, we see that if $n=1$, then $N_1(1,-1)
\in\Om_0(\ga(2/3))$, and $i_1(x,2/3)$ must be even, so
$i_1(x,2/3)=n=1$ is impossible. If $n>1$, we have $n-1>0$
 and $$i_1(x,2/3)=2(k_1+k_2)+n-1.$$ But
$i_1(x,2/3)=n$, so $k_1+k_2=\frac 12$. It is also impossible.

If $k=4$, the solution $(1/2,x)$ itself is a brake orbit. Thus
$i_1(x,1/2)$ and $i_1(x,1)$ are well defined and by Theorem 2.6, we
have that the left hand side equality in (2.8) holds for $k=4$ and
$$i_1(x,1/2)=n,\;\nu_1(x,1/2)=1, i_{L_0}(x,1/4)=i^{L_0}_{\sqrt{-1}}(x,1/4)=0. $$
 By the same arguments as
above, we still get $i_1(x,1/2)=2(k_1+k_2)+n-1$. This is also
impossible.
 \hfill $\blacksquare$\enddemo

\noindent{\bf  Remark 3.2.} If $B=0$, the results of Theorem 3.1
hold  for every $\tau>0$.
  The following  condition is  more accessible than  (HX) but it implies the condition (HX).

\item\;{(H6)} $H''(x)\geq 0$  for all $x\in\r^{2n}$,
the set $D=\{x\in\r^{2n} | H'(x)\not=0,\;0\in\sigma(H''(x))\}$ is
hereditarily disconnected, i.e. every connected component of  $D$
contains only one point.

  Similarly, we consider the brake orbit minimal periodic problem  for the following autonomous
second order Hamiltonian system
$$ \cases & \ddot x+V'(x)=0, \qquad  x\in {\Bbb R}^n,\\
& x(0)=x(\tau/2)=0\\
& x(\tau/2+t)=-x(\tau/2-t), \;\;x(\tau+t)=x(t).\endcases\tag 3.7$$
 A solution $(\tau,x)$ of (3.7) is a kind of brake orbit for the second order Hamiltonian system.

In this paper, we consider  the following conditions on $V$:
\item\quad (V1) $\;\;V\in C^2(\r^n, \r).$
\item\quad (V2) $\;\;$ There exist constants $\mu>2$ and $r_0>0$ such that
$$ 0<\mu V(x)\le V'(x)\cdot x, \qquad  \forall |x|\ge r_0.$$
\item\quad (V3) $\;\; V(x)\ge V(0)=0\;\;\forall x\in \r^n.$
\item\quad (V4) $\;\;V(x)=o(|x|^2), \;\text{at}\; x=0.$
\item\quad (V5) $\;\;$ $V(-x)=V(x),\; \forall x\in \r^n$.
\item\quad (V6) $\;\;$ $V''(x)> 0, \;\;\forall x\in\r. $

\proclaim{Theorem 3.3} Suppose $V$ satisfies the conditions
(V1)-(V6). Then for every $\tau>0$, the problem (3.7) possesses a
non-constant solution $(\tau,x)$ such that the minimal period of $x$
is $\tau$ or $\tau/2$.
\endproclaim

\demo{Proof} Without loss generality, we suppose $\tau=2$. We define
a Hilbert space $W$ which is a subspace of $W^{1,2}([0,1],\r^n)$ by
$$W=\{x\in W^{1,2}([0,1],\r^n)|\,x(t)=\sum_{k=1}^{\infty}\sin k\pi t\cdot a_k, \;a_k\in\r^n\} .$$
The inner product of $W$ is still the $W^{1,2}$ inner product.

 We consider the following functional
$$ \psi (x)=\int^{1}_0 (\frac 1 2 |\dot x|^2-V(x))\,dt,\;\;
                \forall x\in W. \tag 3.8$$
A critical point $x$ of $\psi$ is a solution of the problem (3.7) by
extending  the domain to $\r$ via $x(1+t)=-x(1-t)$ and
$x(2+t)=x(t)$. The condition (V3) implies $\psi(0)=0$. The condition
(V4) implies $\psi(\partial B_{\rho}(0))\ge \alpha_0$ with $\partial
B_{\rho}(0)=\{x\in W\,|\,\|x\|=\rho\}$ for some small $\rho>0$ and
$\alpha_0>0$. In fact, there exists a constant $c_1>0$ such that
$$\int^{1}_0  |\dot x|^2\,dt\ge c_1\|x\|_W^2. \tag 3.9$$
If $\|x\|_W\to 0$, then $\|x\|_{\infty}\to 0$. So by condition (V4),
for any $0<\varepsilon<\frac{c_1}{2}$, there exists small $\rho>0$
such that
$$\int^1_0V(x(t))dt\le \varepsilon\|x\|_2^2\le \varepsilon\|x\|_W^2, \;\;\|x\|_W=\rho.$$
Thus we have
$$\psi (x)=\int^{1}_0 (\frac 1 2 |\dot x|^2-V(x))\,dt\ge (\frac{c_1}{2}-\varepsilon)\rho^2:=\alpha_0>0.$$
The condition (V2) implies that there exists an element $x_0\in W$
with $\|x_0\|>\rho$, such that $\psi(x_0)<0$. In fact, we take an
element $e\in W$ with $\|e\|=1$ and by (V3) we assume
$\int^1_0V(e(t))dt>0$. Consider $x=\lambda e$ for $\lambda>0$.
Condition (V2) implies that there is a constant $c_2>0$ such that
$V(\lambda e)\ge \lambda^{\mu}V(e)-c_2$ for $\lambda$ large enough,
and there holds
$$\psi(\lambda e)\le \lambda^2\int^{1}_0 \frac 1 2 |\dot e|^2dt-\lambda^{\mu}\int^1_0V(e(t))dt+c_2<0.$$
Then we take $x_0=\lambda e$ for large $\lambda$ such that the above
inequalities holds.

 We define
$$\Gamma=\{h\in C([0,1],W)\,|\,h(0)=0,\;\;h(1)=x_0\}$$
and $$c=\inf_{h\in\Gamma}\sup_{s\in [0,1]} \psi(h(s)).$$
 By using the Mountain
pass theorem (cf. Theorem 2.2 of [25]), from the conditions
(V2)-(V4) it is well known that there exists a critical point $x\in
W$ of $\psi$ with critical value $c>0$ which is a Mountain pass
point such that its Morse index satisfying $m^-(x,1)\le 1$. If we
set $y=\dot x$ and $z=(x,y)\in \rr$, the problem (3.7) can be
transformed into the following problem
$$\cases &\dot z=-JH'(z),\\&z(0)\in L_0,\;\;z(1)\in L_0\endcases$$
with $H(z)=H(x,y)=\frac {1}{2} |y|^2+V(x)$. We note that (V5)
implies $H(Nz)=H(z)$, so $(2,z)$ is a brake orbit with brake period
$2$. We remind that in this case the complex structure is $-J$, but
it does not cause any difficult to apply the index theory.
 By Theorem 5.1 of [13], the Morse index $m^-(x,1)$ of $x$ is
just the $L_0$-index $i_{L_0}(z,1)$ of $(1,z)$. i.e., there
holds(see also Lemma 4.6 in the appendix below)
$$m^-(x,1)=i_{L_0}(z,1),\;\;m^0(x,1)=\nu_{L_0}(z,1).$$

We can suppose the minimal period of $x$ is $2/k$ for $k\in\n$.  But
 $i_{L_0}(z,1/k)=m^-(x,1/k)\ge 0$, and from the convexity condition
(V6), we have $i_1(z,2/k)\ge n$. With the same arguments as in the
proof of Theorem 3.1, we get $k\in\{1,2\}$. \hfill
$\blacksquare$\enddemo

We note that the functional $\psi$ is even, there may be infinite
many solutions $(\tau,x)$ satisfying  Theorem 3.3. We also note that
Theorem 3.3 is not  a special case of Theorem 3.1, since the
Hamiltonian function $H(x,y)=\frac 12|y|^2+V(x)$ is quadratic in the
variables $y$, in this case $B=\pmatrix 0 & 0\\0 & I_n\endpmatrix$
with $\|B\|=1$. Thus when applying Theorem 3.1 to this case, we can
only get the result of Theorem 3.3 for $0<\tau<2\pi$.

We now consider the following problem

$$ \cases & \ddot x+V'(x)=0, \qquad  x\in {\Bbb R}^n,\\
& \dot x(0)=\dot x(\tau/2)=0,\\
& x(\tau/2+t)=x(\tau/2-t), \;\;x(\tau+t)=x(t).\endcases\tag 3.10$$
 A solution of (3.10) is also a kind of brake orbit for the second order Hamiltonian system.

By set $y=\dot x$, $z=(y,x)$ and $H(z)=H(y,x)=\frac 12|y|^2+V(x)$,
the problem (3.10) can be transformed into the following
$L_0$-boundary value problem
$$\cases\dot z=JH'(z)\\
z(0)\in L_0,\;\;z(\tau/2)\in L_0.\endcases$$
 In this case the condition $H(Nz)=H(z)$ is satisfied automatically.
 Set $B=\pmatrix I_n &0\\0& 0\endpmatrix$, then $\|B\|=1$. The
 following result is a direct consequence of Theorem 3.1.

\proclaim{Corollary 3.4} Suppose $V$ satisfies the conditions
(V1)-(V4) and (V6). Then for every $0<\tau<2\pi$, the problem (3.10)
possesses a non-constant  solution $(\tau,x)$ such that $x$ has
minimal period $\tau$ or $\tau/2$.
\endproclaim

We note that if we directly solve the problem (3.10) by the same way
as in the proof of Theorem 3.3, the formation of the functional is
still $\psi$ as defined in (3.8), but the domain should be
$$W_1=\{x\in W^{1,2}([0,1],\r^n)|\,x(t)=\sum_{k=0}^{\infty}\cos k\pi t\cdot a_k, \;a_k\in\r^n\}.$$
In this time, it is not able to apply the Mountain pass theorem to
get a critical point  directly due to the fact $\r^n\subset W_1$, so
the inequality (3.9) is not true.

\head \S4 Appendix. Some properties for the indeices \endhead
 \subhead 4.1. Some properties of Maslov-type index \endsubhead For a
 symplectic path $\gamma\in \Cal {P}(2n)$, its Maslov-type index
  is a pair of integers $(i_1(\gamma),\nu_1(\gamma))\in \Bbb{Z}\times
  \{0,1,\cdots,2n\}$ (cf. [20],[21]). If $\gamma\in \Cal {P}(2n)$ is the fundamental
  solution of a linear Hamiltonian system
  $$\dot x=JB(t)x$$
  with  continuous symmetric matrix  function $B(t)$, its
  Maslov-type index usually denoted by $(i_1(B),\nu_1(B))$. The
  following result was proved in [12].
\proclaim{Lemma 4.1} If $B_1(t)-B_2(t)>0$ is a positive definite
matrix function, then there holds
$$i_1(B_1)\ge i_1(B_2)+\nu_1(B_2). \tag 4.1$$
(4.1) also holds under the following condition
$$B(t)=B_1(t)-B_2(t)\ge 0,\;\;\;\int_0^1B(t)dt>0.$$\endproclaim
As a direct consequence, if the continuous symmetric matrix function
satisfying $B(t)\ge 0$ and $\int^1_0B(t)dt>0$, then there holds
$$i_1(B)\ge n. \tag 4.2$$
\proclaim{Definition 4.2}([20],[21])Two symplectic paths $\gamma_0$
and $\gamma_1\in \Cal P(2n)$ are
 homotopic on $[0,1]$, denoted by $\gamma_0\sim \gamma_1$, if there exists a map $\delta\in C([0,1]\times [0,1], \Sp(2n))$ such that
 $\delta(0,\cdot)=\gamma_0(\cdot),\;\delta(1,\cdot)=\gamma_1(\cdot),\;\delta(s,\cdot)=I_{2n}$,
 and $\nu_1(\delta(s,1))$ is constant for $0\le s\le 1$.  \endproclaim

We note that for two paths $\gamma_0$ and $\gamma_1\in \Cal P(2n)$
with the same end points $\gamma_0(1)=\gamma_1(1)$, $\gamma_0\sim
\gamma_1$ with fixed end points if and only if
$i_1(\gamma_0)=i_1(\gamma_1)$. By choosing suitable {\it zigzag
standard paths} $\alpha_{n,k}$ in $\Cal P^*(2n)$ with
$i_1(\alpha_{n,k})=k$ and $\alpha_{n,k}(1)=M_n^{\pm}$ if $(-1)^k=\pm
1$  as in [22], and by the definition of the Maslov-type index, we
have the following result.
 \proclaim{Lemma 4.3} For a symplectic path $\gamma\in \Cal P(2n)$ with $\gamma(1)=M_1\diamond M_2$,
 $M_j\in \Sp(2n_j)$, $j=1,2$, $n_1+n_2=n$, there exists two symplectic paths
 $\gamma_j\in \Cal P(2n_j)$ such that $\gamma\sim \gamma_1\diamond
 \gamma_2$ and $\gamma_j(1)=M_j$.
 \endproclaim
The index function $(i_{\omega}(\gamma),\nu_{\omega}(\gamma))$ was
defined for $\omega\in \Bbb U:=\{z\in \Bbb C|\;|z|=1\}$ in [22] by
Y.Long.

\proclaim{Definition 4.4}([22])For any $M\in \Sp(2n)$ and $\omega\in
\Bbb U$, choosing $\gamma\in \Cal P(2n)$ with $\gamma(1)=M$, the
splitting numbers of $M$ are defined by
$$S^{\pm}_M(\omega)=\lim_{\epsilon\to 0^+}i_{\exp(\pm \epsilon\sqrt{-1})\omega}(\gamma)-i_{\omega}(\gamma).$$\endproclaim
The following list for the splitting number comes from [21].

\item {(1)} $(S^+_M(1),S^-_M(1))=(1,1)$ for $M=N_1(1,b)$ with $b=1$ or
$0$.

\item {(2)} $(S^+_M(1),S^-_M(1))=(0,0)$ for $M=N_1(1,-1)$.

\item {(3)} $(S^+_M(-1),S^-_M(-1))=(1,1)$ for $M=N_1(-1,b)$ with $b=-1$ or
$0$.

\item {(4)} $(S^+_M(-1),S^-_M(-1))=(0,0)$ for $M=N_1(-1,1)$.

\item {(5)} $(S^+_M(e^{\sqrt{-1}\theta}),S^-_M(e^{\sqrt{-1}\theta}))=(0,1)$ for $M=R(\theta)$
with $\theta\in (0,\pi)\cup (\pi,2\pi)$.

\item {(6)} $(S^+_M(\omega),S^-_M(\omega))=(1,1)$ for $M=N_2(\omega,b)$ being non-trivial with $\omega=e^{\sqrt{-1}\theta}\in \Bbb U
\setminus\Bbb R$.

\item {(7)} $(S^+_M(\omega),S^-_M(\omega))=(0,0)$ for $M=N_2(\omega,b)$ being trivial with $\omega=e^{\sqrt{-1}\theta}\in \Bbb U
\setminus\Bbb R$.

\item {(8)} $(S^+_M(\omega),S^-_M(\omega))=(0,0)$ for any $\omega\in \Bbb U$ and $M\in \Sp(2n)$ satisfying
$\sigma(M)\cap \Bbb U=\emptyset$.

\subhead 4.2. Some properties of the $L_0$-index\endsubhead For a
symplectic path $\gamma\in \Cal P(2n)$, the so called $L_0$-index
$(i_{L_0}(\gamma),\nu_{L_0}(\gamma))\in \Bbb Z\times
\{0,1,\cdots,n\}$ was  first defined in [13]. We have a brief
introduction of this index theory in the section 2 of this paper.
The following result was proved in [13].
 \proclaim{Lemma 4.5} Suppose $\gamma\in \Cal P(2n)$ is the
 fundamental solution of the following linear Hamiltonian system
 $$\dot x(t)=JB(t)x(t),\;\;x(t)\in \Bbb R^{2n},$$
 where $B(t)=\pmatrix S_{11}(t) & S_{12}(t)\\S_{21}(t) &
 S_{22}(t)\endpmatrix$ is symmetric with $n\times n$ blocks
 $S_{jk}$. If $S_{22}(t)>0$(positive definite), there holds
 $$i_{L_0}(\gamma)=i_{L_0}(B)\ge 0. \tag 4.3$$
 (4.3) is also true if $S_{22}(t)\ge 0$ and $\int^1_0S_{22}(t)dt>0$.
\endproclaim

We consider the following problem
$$\cases & [P(t) x'(t)-Q(t)x(t)]'+Q^T(t)x'(t)+R(t)x=0,\\& x(0)=x(1)=0,\endcases\tag 4.4 $$
where $P$ and $R$ are symmetrial $n\times n$ matrix function, we
suppose $-P>0$ (positive definite). For simplicity, We assume $P, Q
$ are smooth  and $R$ is continuous. The equations in (4.4) was
studied by M.Morse.  We turn it into a first order equations with
Lagrangian boundary condition by setting $z(t)=(x(t), y(t))^T\in\rr$
with $y=P(t)x'(t)-Q(t)x(t)$:
$$\cases & \dot z=JB(t)z\\& z(0), \;z(1)\in L_0, \endcases\tag 4.5 $$
where $B=B(t)$ is defined by
$$B(t)=\pmatrix -R(t)-Q^T(t)P^{-1}(t)Q(t) & -Q^T(t)P^{-1}(t)\\-P^{-1}(t)Q(t) & -P^{-1}(t)\endpmatrix.   $$
We take the space $W=W_0^{1,2}([0,1], \r^n)$, the subspace of
$W^{1,2}([0,1], \r^n)$ with the elements $x$ satisfying
$x(0)=x(1)=0$. Define the following functional on $W$
$$\aligned\varphi(x)&=-\frac12\int_0^1 \<P^{-1}(t)(P(t)x'(t)-Q(t)x(t)), P(t)x'(t)-Q(t)x(t)\>\\&
-\<(R(t)+Q^T(t)P^{-1}(t)Q(t))x(t),x(t)\>\,dt.\endaligned $$ The
critical point of $\varphi$ is a solution of the problem (4.4), and
so we get a solution of the problem (4.5).
 Denote the Morse index of the functional $\varphi$
at $x=0$ by $m^{L_0}(B)$, which is the total multiplicity of the
negative eigenvalues of the Hessian of $\varphi$ at $x=0$, and the
nullity by $n^{L_0}(B)$. The following result was proved in [13].
 \proclaim{Lemma 4.6} There holds
 $$i_{L_0}(B)=m^{L_0}(B),\;\;  \nu_{L_0}(B)=n^{L_0}(B). $$
 \endproclaim

 Let $E$ be a separable Hilbert space, and $Q=A-B: E\to E$ be a
bounded salf-adjoint linear operators with $B:E\to E$ a compact
self-adjoint operator. $N=\ker Q$ and $\dim N<+\infty$.
$Q|_{N^{\bot}}$ is invertible. $P:E\to N$ the orthogonal projection.
Set $d=\frac 14 \|(Q|_{N^{\bot}})^{-1}\|^{-1}$.
$\Gamma=\{P_k|k=1,2,\cdots\}$ be the Galerkin approximation sequence
of $A$:

\item\;(1)  $E_k:=P_kE$ is finite dimensional for all $k\in\n$,

\item\;(2)  $P_k\to I$ strongly as $k\to +\infty$

\item\;(3)  $P_kA=AP_k$.

For an operator $S$, we denote by $M^{*}(S)$ the eigenspaces of $S$
with eigenvalues belonging to $(0,+\infty)$, $\{0\}$ and $(-\infty,
0)$ with  $*=+,0$ and $*=-$, respectively. We denote by $m^*(S)=\dim
M^*(S)$. Similarly, we denote by $M_d^{*}(S)$ the $d$-eigenspaces of
$S$ with eigenvalues belonging to $(d,+\infty)$, $(-d,d)$ and
$(-\infty, -d)$ with  $*=+,0$ and $*=-$, respectively. We denote by
$m_d^*(S)=\dim M_d^*(S)$.

 \proclaim{Lemma 4.7}(Lemma 3.3 in [16]) Let $B$ be a linear symmetric  compact
 operator. Then the difference of the $d$-Morse
 indices
 $$m_d^-(P_m(A-B)P_m)-m_d^-(P_mAP_m) \tag 4.6$$
 eventually becomes a constant independent of $m$, where $d>0$ is determined by the operators $A$ and $A-B$.
 Moreover $m^0_d(P_m(A-B)P_m)$ eventually becomes a constant independent of
 $m$ and for large $m$, there holds
$$m_d^0(P_m(A-B)P_m)=m^0(A-B). \tag 4.7$$
\endproclaim
\proclaim{Definition 4.8}([16]) For the operators $A$ and $B$ in
Lemma 4.7, $\Gamma$ is an Galerkin approximation sequence w.r.t $A$,
we define the relative index by
$$I(A,A-B)= m_d^-(P_m(A-B)P_m)-m_d^-(P_mAP_m), m\ge m^*, \tag 4.8$$
where $m^*>0$ is a constant large enough such that the  difference
in (4.6) becomes a constant independent of $m\ge m^*$.
\endproclaim
For $\omega=e^{\sqrt{-1}\theta}$, we define a Hilbert space
$E^{\omega}=E^{\omega}_{L_0}$ consisting of those $x(t)$ in
$L^2([0,1], \c^{2n})$ such that $e^{-\theta t J}x(t)$ has Fourier
series
$$e^{-\theta t J}x(t)=\sum_{j\in \z}e^{j\pi tJ}\pmatrix 0\\a_j\endpmatrix,\;a_j\in \c^n  $$
and
$$\|x\|^2:=\sum_{k\in\z}(1+|k|)|a_k|^2<\infty. $$
For $x\in E^{\omega}$, we can write
$$\align x(t)&=e^{\theta tJ}\sum_{j\in\z}e^{j\pi tJ}\pmatrix 0\\a_j\endpmatrix
 =\sum_{j\in\z}e^{(\theta+j\pi)tJ}\pmatrix 0\\a_j\endpmatrix \\
 &=\sum_{j\in\z}e^{(\theta+j\pi)t\sqrt{-1}}\pmatrix \sqrt{-1}a_j/2\\a_j/2\endpmatrix+
 e^{-(\theta+j\pi)t\sqrt{-1}}\pmatrix -\sqrt{-1}a_j/2\\a_j/2\endpmatrix.\endalign $$
So we can write
$$x(t)=\xi(t)+N\xi(-t), \; \xi(t)=\sum_{j\in\z}e^{(\theta+j\pi)t\sqrt{-1}}
\pmatrix \sqrt{-1}a_j/2\\a_j/2\endpmatrix. \tag 4.9$$ For
$\omega=1$, i.e., $\theta=0$, we define two selfadjoint operators
$A_1, B_1\in \Cal L(E^1)$ by extending the bilinear forms
$$\<A_1x,y\>=\int^1_0(-J\dot x(t),y(t))dt,\;\;\<B_1x,y\>=\int^1_0(B(t)x,y)dt $$
on $E^1$, here $(\cdot,\cdot)$ is the Hermitian inner product in
$\c^{2n}$. Then $B$ is compact. For $\omega=e^{\sqrt{-1}\theta}$,
$\theta\in [0,\pi)$, we define two self-adjoint operators
$A^{\omega}, B^{\omega}\in \Cal L(E^{\omega})$ by extending the
bilinear forms
$$\aligned & \<A^{\omega}x,y\>=\int^1_0(-J\dot x(t),y(t))dt,\;\;\\&\<B^{\omega}x,y\>=\int^1_0(B(t)x(t),y(t))dt \endaligned$$
 on $E^{\omega}$, where we have written $x(t)=\xi(t)+N\xi(-t),\; y(t)=\eta(t)+N\eta(-t)$ as in (3.10). Then $B^{\omega}$ is also compact.

 By Theorem 2.1 of [14], we have the following formula
 $$I(A_1,A_1-B_1)=i_{L_0}(B)+n. \tag 4.10$$
 \proclaim{Definition 4.9}([16]) We define the index function
 $$i_{\omega}^{L_0}(B):=I(A^{\omega}, A^{\omega}-B^{\omega}),\;\;\nu_{\omega}^{L_0}(B):=m^0(A^{\omega}-B^{\omega}),\;
  \omega=e^{\sqrt{-1}\theta},\;\theta\in (0,\pi). $$
 \endproclaim

 \proclaim{Lemma 4.10}([16]) The index function $i_{\omega}^{L_0}(B)$ is
locally constant. For $\omega_0=e^{\sqrt{-1}\theta_0},\;\theta_0\in
(0,\pi)$ is a point of discontinuity of $i_{\omega}^{L_0}(B)$, then
$\nu_{\omega_0}^{L_0}(B)>0$ and so $\dim (\gamma(1)L_0\cap
e^{\theta_0 J}L_0)>0$. Moreover there hold
$$\aligned & |i_{\omega_0+}^{L_0}(B)-i_{\omega_0-}^{L_0}(B)|\le \nu_{\omega_0}^{L_0}(B),\;
|i_{\omega_0+}^{L_0}(B)-i_{\omega_0}^{L_0}(B)|\le
\nu_{\omega_0}^{L_0}(B),\\&
|i_{\omega_0-}^{L_0}(B)-i_{\omega_0}^{L_0}(B)|\le
\nu_{\omega_0}^{L_0}(B), |i_{L_0}(B)+n-i_{0+}^{L_0}(B)|\le
\nu_{L_0}(B), \endaligned \tag 4.11$$
 where $i_{\omega_0+}^{L_0}(B)$, $i_{\omega_0-}^{L_0}(B)$ are the
 right and left limit respectively of the index function $i_{\omega}^{L_0}(B)$
 at $\omega_0=e^{\sqrt{-1}\theta_0}$ as a function of $\theta$.
\endproclaim
By (4.10), Definition 4.9 and Lemma 4.10, we see that for any
$\omega_0=e^{\sqrt{-1}\theta_0},\;\theta_0\in (0,\pi)$, there holds
$$i^{L_0}_{\omega_0}(B)\ge i_{L_0}(B)+n-\sum\Sb\omega=e^{\sqrt{-1}\theta}\\
0\le \theta\le \theta_0\endSb\nu^{L_0}_{\omega}(B). \tag 4.12 $$
 We note that
 $$\sum\Sb\omega=e^{\sqrt{-1}\theta}\\
0\le \theta\le \theta_0\endSb\nu^{L_0}_{\omega}(B)\le n. \tag 4.13
$$
So we have $$i_{L_0}(B)\le i^{L_0}_{\omega_0}(B)\le i_{L_0}(B)+n.
\tag 4.14$$

 {\bf Acknowledgements}: The author  of this paper appreciates the
referee for his valuable suggestions.


\Refs \widestnumber\key{ASD1}

\ref\key 1\by V. Benci and F. Giannoni\paper A new proof of the
existence of a brake orbit. In ``Advanced Topics in the Theory of
Dynamical Systems". \jour Notes Rep. Math. Sci. Eng.\vol 6 \yr
1989\pages 37-49\endref

\ref\key 2 \by  K. C. Chang \book  Infinite Dimensional Morse Theory
and Multiple Solution Problems \publ Birkh\"auser \publaddr Basel
\yr 1993
\endref



\ref\key  3 \by     C. Conley and E. Zehnder \paper Maslov-type
index theory for flows and periodic solutions for Hamiltonian
equations \jour  Commun. Pure Appl. Math. \vol 37 \yr    1984 \pages
207-253
\endref

\ref\key 4 \by D. Dong and Y. Long \paper The Iteration Formula of
the Maslov-type Index Theory with Applications to Nonlinear
Hamiltonian Systems \jour Trans. Amer. Math. Soc. \vol  349  \yr
1997   \pages 2619-2661
\endref

\ref\key 5 \by I. Ekeland \book Convexity Methods in Hamiltonian
Mechanics \publ Springer \publaddr Berlin \yr 1990
\endref

\ref\key 6 \by I. Ekeland and H. Hofer \paper Periodic solutions
with prescribe period for convex autonomous Hamiltonian systems
\jour Invent. Math. \vol 81 \yr 1985\pages 155-188\endref

\ref\key 7
 \by G. Fei
 \paper Relative Morse index and its application to Hamiltonian
 systems in the Presence of symmetries\jour J.Diff.Equ. \vol 122\yr
 1995\pages 302-315\endref

\ref\key 8 \by G. Fei and Q. Qiu \paper Minimal period solutions of
nonlinear Hamiltonian systems \jour Preprint 1996
\endref

\ref\key 9
 \by G. Fei and Q. Qiu
 \paper Periodic solutions of asymptotically linear Hamiltonian
 systems
 \jour Chin. Ann. of Math. \vol 18B:3 \yr 1997\pages 359-372\endref



\ref\key 10 \by M. Girardi and M. Matzeu \paper Periodic solutions
of convex Hamiltonian systems with a quadratic growth at the origin
and superquadratic at infinity \jour Ann. Math. Pura ed App. \vol
147 \yr 1987 \pages 21-72
\endref

\ref\key 11 \by C. Li and C. Liu
 \paper Nontrivial solutions of
superquadratic Hamiltonian systems with Lagrangian boundary
conditions and the L-index theory\jour Chinese Ann.Math. \vol
29(6)\yr 2008\pages  597-610\endref

\ref\key 12
  \by C. Liu
  \paper  A note on the monotonicity of Maslov-type index of
linear Hamiltonian systems with applications\jour Proceedings of the
Royal Society of Edinburgh\vol 135A\pages 1263-1277\yr 2005\endref

\ref\key 13
  \by C. Liu
  \paper Maslov-type index theory for symplectic
paths with Lagrangian boundary conditions \jour Advanced Nonlinear
Studies \vol 7 \yr 2007 \pages 131-161\endref

 \ref\key 14
  \by C. Liu
  \paper  Asymptotically linear Hamiltonian systems with Lagrangian
boundary conditions \jour Pacific J. Math(in press)
  \endref



\ref\key 15
 \by C. Liu and Y. Long
 \paper Iteration inequalities of the Maslov-type index theory with applications
 \jour J.Diff.Equa. \vol 165\yr 2000\pages 355-376\endref

 \ref\key 16
 \by C. Liu and D. Zhang
 \paper An iteration theory of Maslov-type index for symplectic paths associated with a Lagrangian
subspace  and Multiplicity of brake orbits in bounded convex
symmetric domains \jour arXiv:0908.0021v1 [math.SG]\endref

\ref\key 17 \by Y. Long \paper The minimal period problem of
classical Hamiltonian systems with even potentials \jour Ann. I. H.
P. Anal. non lin\'eaire \vol 10-6 \yr 1993 \pages 605-626
\endref

\ref\key 18 \by Y. Long \paper The minimal period problem of
periodic solutions for autonomous superquadratic second order
Hamiltonian systems \jour J. Diff. Equa. \vol 111 \yr 1994 \pages
147-174
\endref

\ref\key 19 \by Y. Long \paper  On the minimal period for periodic
solution problem of nonlinear Hamiltonian systems \jour Chinese Ann.
of  Math. \vol 18B \yr 1997 \pages 481-484
\endref

\ref\key 20 \by Y. Long \paper Maslov-type index, degenerate
critical points, and asymptotically linear hamiltonian systems \jour
Science in China(Series A) \vol 33 \yr 1990 \pages 1409-1419\endref

 \ref\key 21
 \by Y. Long
 \book Index Theory for Symplectic Path with Applications
  \publ Birkh\"auser Verlag, Basel, Boston, Berlin \yr 2002
 \endref

 \ref\key 22
 \by Y. Long
 \paper Bott formula of the Maslov-type index theory
 \jour Pacific J. Math \yr 1999\vol 187\pages 113-149\endref



\ref\key 23
 \by Y. Long, D. Zhang and C. Zhu
 \paper Multiple brake orbits in bounded convex symmetric domains
 \jour Advances in Mathematics \vol 203 \pages 568-635\yr 2006
 \endref



\ref\key 24 \by P. Rabinowitz \paper Periodic solutions of
Hamiltonian systems \jour Comm. Pure Appl. Math. \vol 31 \yr 1978
\pages 157-184
\endref

\ref\key 25 \by P. Rabinowitz \paper Minimax methods in critical
point theory with applications to differential equations \jour CBMS
Regional Conf. Ser. in Math. Amer. Math. Soc. \vol 65  \yr 1986
\endref


\ref\key 26 \by  P.  Rabinowitz
 \paper  On the existence of periodic solutions for a
class of symmetric Hamiltonian systems\jour Nonlinear Anal. T. M.
A.\vol 11 \yr 1987\pages 599-611\endref


\ref\key 27\by  A. Szulkin\paper An index theory and existence of
multiple brake orbits for star-shaped Hamiltonian systems\jour Math.
Ann.\vol 283 \yr 1989\pages 241-255\endref


\endRefs
\enddocument